\documentclass[11pt,letterpaper]{amsart}

\usepackage[utf8]{inputenc}
\usepackage{subcaption}
\usepackage{verbatim} 

\usepackage{graphicx} 

\usepackage{amsmath}
\usepackage{amsfonts}
\usepackage{amssymb}
\usepackage{hyperref}
\newtheorem{theorem}{Theorem}
\newcommand{\R}{\mathbb{R}}
\newcommand{\supp}{\textnormal{supp }}

\usepackage[backend=biber,style=numeric]{biblatex}
\title[TAT with Circular Integrating Detectors]{Thermoacoustic Tomography with Circular Integrating Detectors and Variable Wave Speed}
\author[C. Mathison]{Chase Mathison}

\address{Department of Mathematics, Purdue University, West Lafayette, IN 47907}

\thanks{Author partly supported by NSF Grant DMS-1600327}

\date{\today} 

\begin{document}
\begin{abstract}
    We explore Thermoacoustic Tomography with circular integrating detectors assuming variable, smooth wave speed.  We show that the measurement operator in this case is a Fourier Integral Operator and examine how the singularities in initial data and measured data are related through the canonical relation of this operator.  We prove which of those singularities in the initial data are visible from a fixed open subset of the set on which measurements are taken.  In addition, numerical results are shown for both full and partial data.
\end{abstract}
\maketitle
\section{Introduction}

Thermoacoustic Tomography is a medical imaging method in which a short pulse of electromagnetic radiation is used to excite cells in some object we wish to image, typically the organs of a patient.  Upon absorbing the EM radiation, the cells in the patient in turn vibrate, creating ultrasonic waves that then propagate out of the patient and are measured by any number of methods.  Using this measured data, we then try to reconstruct, in some sense, an image of the inside of the patient.  This is a hybrid imaging method which uses high contrast, low resolution EM radiation to excite the cells, and low contrast, high resolution ultrasound waves as measurement \cite{ora}.  The hope is to be able to get an image with good contrast and resolution by combining these two types of waves.  The case of point-wise measurements with constant and variable wave speed in the region of interest has been studied extensively \cite{kuchment,stefAndUhl,xu}.  Other methods of measurement of the ultrasonic waves include measurements with linear integrating detectors \cite{kruger1}, planar integrating detectors \cite{haltmeier, stefanov} and circular integrating detectors or cylindrical stacks of circular integrating detectors \cite{zangerl2,zangerl}.  Circular integrating detectors have a few advantages over linear integrating detectors and planar integrating detectors, including compactness of the experimental setup \cite{zangerl2}.  The case of planar integrating detectors was studied in \cite{stefanov}, and that work focused on the problem with a smooth, variable wave speed.  The case of circular (and cylindrical) integrating detectors with constant wave speed has been studied in \cite{zangerl2,zangerl}.  In those works, explicit formulae are given for reconstruction of an initial pressure density using full measurements, i.e. measurements for every circular integrating detector of a fixed radius with center on the unit circle, for all time.  That reconstruction is stable in the case that the object being imaged is contained in the interior of the circular integrating detectors, but is unstable for the case that the object lies entirely outside of the circular integrating detector.  The present work focuses on the case of circular integrating detectors in the plane with a 2 dimensional region of interest.  Further, we do not make a constant wave speed assumption, we only assume that the wave speed $c(x) > 0$ is smooth in all of $\R^2$ and is $1$ outside of a compact subset of $\R^2$.  We show that the measurement operator in this case is a Fourier Integral Operator and compute its canonical relation, which happens to be a local dif{}feomorphism, thus allowing us to determine how singularities in initial data propagate to the measurement data.  We also show that this operator is injective and prove stability of the measurement operator, and in addition we prove what singularities in the initial data are visible from a fixed open subset of the set of points on the circle where the measurements are taken in a given time interval.  Lastly, we provide numerical results obtained through simulation in Matlab using both full and partial data that support our findings.

\section{Setup}
We begin by defining the space of distributions that our initial pressure distribution must be in.  Let $$H_D(\Omega) = \left\{ f \in \mathcal{E}'(\Omega) \mid \int_{\Omega} |\nabla f|^2\,dx < \infty \right\},$$ where $\Omega \subset \R^2$ is open and $\mathcal{E}'(\Omega)$ is the space of distributions compactly supported in $\Omega$.  This is the natural space in which to take $f$ when the energy of the system is taken into consideration. Let $\lVert f \rVert_{H_D(\Omega)}^2 = \int_\Omega |\nabla f|^2 \,dx$.  The space $H_D(\Omega)$ is the completion of $C_0^\infty (\Omega)$ under the given norm. We further suppose that $\bar{\Omega} \subset B_1(0)$ where $B_1(0)$ is the unit ball in $\mathbb{R}^2$ centered at the origin.  We view $f$ as an initial pressure distribution of some object to be imaged represented by $\Omega$.  Then, after exposing $\Omega$ to microwave radiation, the ultrasonic waves created solve the acoustic wave equation given by
\begin{equation}\label{eq:1}
\begin{cases}
\partial_t^2 u(t,x) = c^2(x)\Delta u(t,x), & (t,x) \in \mathbb{R}\times \mathbb{R}^2 \\
u(0,x) = f(x), & f(x) \in H_D\left( \Omega \right) \\
\partial_t u(t,x) \mid_{t=0} = 0, & x \in \R^2
\end{cases}
\end{equation}
where $0 < c(x) \in C^\infty(\R^2)$ is the smooth wave speed, assumed to be known.  Outside of $\Omega$, $c \equiv 1$.  The problem of interest is to detect these waves, solutions $u(t,x)$ to the above wave equation, with detectors located on the boundary of $B_R(0)$, and then using these measurements, reconstruct the initial pressure distribution $f$.  As mentioned in the introduction, extensive research has been done in the constant speed case ($c(x) \equiv 1$ for all $x \in \mathbb{R}^2$) and variable speed case with point detectors in which we assume access to $u\mid_{U\times \Gamma}$ where $\Gamma \subset S^1$ is open and $U \subset \R$ is some time interval.  Research has also been done for linear and planar integrating detectors in both the constant and variable speed case, and also circular and cylindrical integrating detectors in the constant wave speed case.  When imaging with these integrating detectors, instead of assuming direct access to $u$ on some open subset of the boundary of $\Omega$, the measured data is an average of $u$ over a circular detector of radius $r$ centered on the boundary of the ball of radius $R$ (which we will choose later), and data is assumed to be collected on an open subset of of this boundary, not necessarily the entirety of the boundary.  The present work considers the problem with variable speed in $\Omega$, and circular integrating detectors.  We will have two cases to consider, which we will call the large radius detector case and the small radius detector case, which are depicted in Figure \ref{fig:cases}.  The large radius detector case is the experimental setup in which $\Omega$ is on the ``inside'' of the circular integrating detectors, and the small radius detector case is the setup in which $\Omega$ is on the ``outside'' of the circular integrating detector.

\section{Construction}
We are interested in seeing what singularities we can recover from rotating a circular integrating detector around some object that has been heated via microwaves.  To start, we recall that solving the wave equation (as above) up to a smooth error, for $x \in \mathbb{R}^2$, can be accomplished by use of the geometric optics construction (see Section 3 of \cite{taylor81}).  After making this construction, we have that $u$ up to a smooth error is given locally by

$$
u(t,x) = \frac{1}{(2\pi)^2} \sum_{\sigma = \pm} \int e^{i\phi_\sigma (t,x,\xi)}a_\sigma(t,x,\xi) \hat{f}(\xi)\,d\xi.
$$
Here, $\phi_\sigma$ solves the \emph{eikonal equation}: $((\phi_\sigma)_t)^2 = c^2(x)|\nabla_x \phi_\sigma|^2$ with initial condition $\phi_\sigma(0,x,\xi) = x\cdot\xi$.  In addition, $\phi_\sigma$ is positive homogeneous of order $1$ in $\xi$, i.e. $\phi_\sigma(t,x,\lambda \xi) = \lambda \phi_\sigma (t,x,\xi)$ for $\lambda > 0$.  The eikonal equation is only solvable locally in time, which results in our solution $u(t,x)$ being only a local solution in time.  This is not actually an issue however, as we simply repeat this procedure with new initial conditions.  The resulting Fourier Integral Operator is then actually a composition of Fourier Integral Operators, which is a technical issue we refer the reader to \cite{stefAndUhl} for more details.  Because of this, we may assume that the eikonal equation is solvable until geodesics intersect circular integrating detectors.  By this construction, we may obtain $a_\sigma$ as a classical symbol of order $0$ with expansion $a_\sigma(t,x,\xi) \sim \sum_{j\ge 0} a_{\sigma}^{(j)}(t,x,\xi)$ where $a_{\sigma}^{(j)}(t,x,\xi)$ is positively homogeneous of order $-j$ for $|\xi|$ large.  In particular, $a_{\sigma}^{(0)}(t,x,\xi)$ solves the transport equation
$$
\left( (\partial_t \phi_\sigma)\partial_t - c^2(x) (\nabla_x \phi_\sigma)\cdot \nabla_x + \frac{1}{2}(\partial_t^2 - c^2(x)\Delta_x ) \phi_\sigma\right) a_\sigma^{(0)} = 0
$$
with initial conditions $a_\sigma^{(0)}(0,x,\xi) = 1/2$.  The last term on the left hand side of this equation acts on $a_\sigma^{(0)}$ by multiplication.

\begin{figure}[h!]
    \centering
    \includegraphics{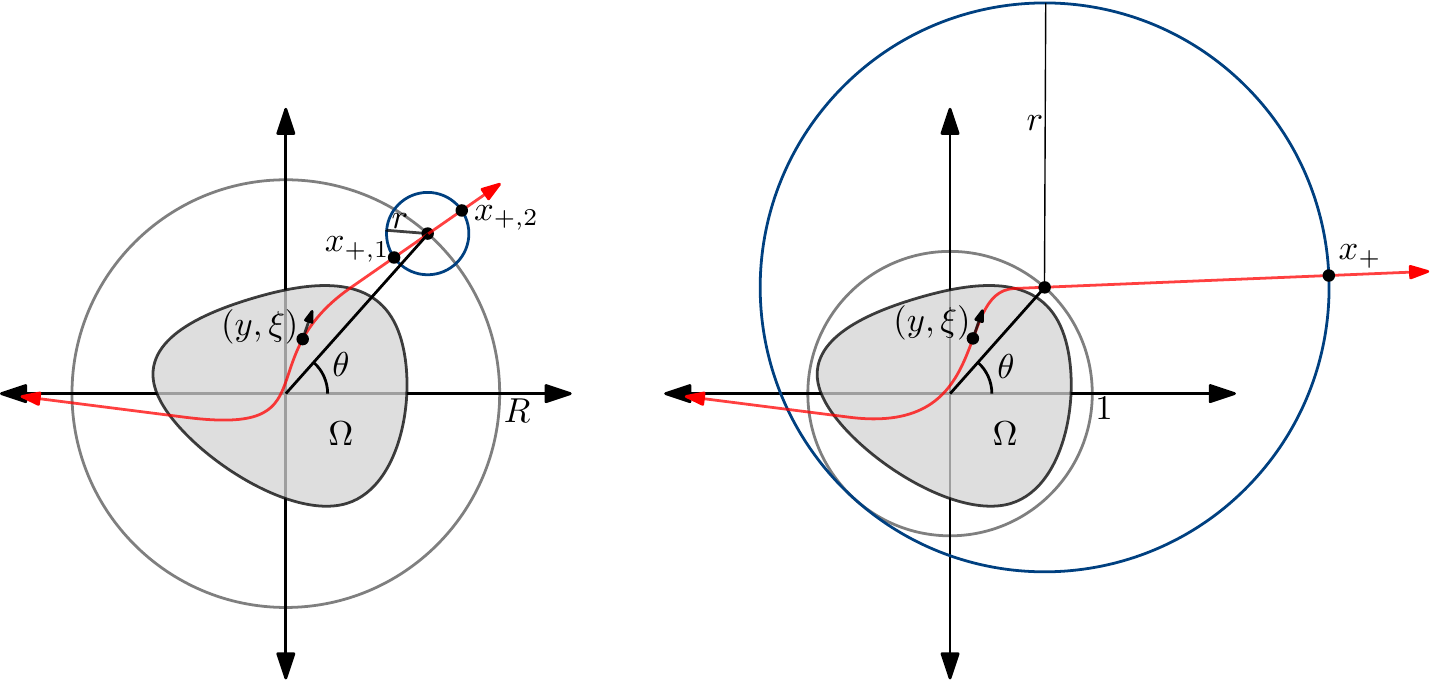}
    \caption{Two dif{}ferent experimental setups shown depending on the radius of the integrating detector. On the left is the small radius case, and on the right is the large radius case.}
    \label{fig:cases}
\end{figure}

Now in the situation of Thermoacoustic Tomography using circular integrating detectors around the object we wish to image, the measurements at the detector are given by the circular Radon transform:

$$
Mf(t,\theta) = \frac{1}{2\pi} \int_0^{2\pi} u(t,C(\theta, \alpha))\,d\alpha,
$$
where $C(\theta, \alpha) = R(\cos\theta,\sin\theta) + r(\cos\alpha,\sin\alpha)$ is a parametrization of the circular detector, $R$ is the distance from the origin to the center of the circle, $r$ is the radius of the circular detector, $\theta$ is the angle made between the positive horizontal axis and the ray from the origin to the center of the circle, and $u(t,x)$ is the solution to the IVP (\ref{eq:1}).  The radius of the circular integrating detector, $r$, must be chosen so that the detector does not intersect $\Omega$.  To accomplish this, we must have either $r$ small enough so that $R - r \ge 1$, guaranteeing the detector does not intersect $\Omega \subset B_1(0)$, or, we could fix $R = 1$ and choose $r \ge 2$, in which case $\Omega$ is contained in the interior of the disc defined by the detector (see Figure \ref{fig:cases}).  For convenience, we define $\boldsymbol{\theta} = ( \cos \theta,  \sin \theta)$.  We can rewrite $Mf$ by using the $\delta$ distribution:

$$
Mf(t,\theta) = \frac{1}{2\pi}\int_{\R^2} u(t,x)\delta(\left| x - R\boldsymbol{\theta} \right|^2 - r^2)\,dx.
$$

We now plug in our solution for $u$ obtained via the geometric optics construction and denote by $M_+$ and $M_-$ the operators taking $\sigma = +$ and $\sigma = -$ respectively after substituting the geometric optics solution in $Mf(t,\theta)$.  Then $Mf = M_+f + M_-f$, where
\begin{gather*}
    M_+ f = \frac{1}{(2\pi)^3} \iint_{\R^2\times \R^2} e^{i\phi_{+}(t,x,\xi)}a_+(t,x,\xi)\hat{f}(\xi)\,d\xi\delta(|x-R\boldsymbol{\theta}|^2-r^2)\,dx \\
    M_- f = \frac{1}{(2\pi)^3} \iint_{\R^2\times \R^2} e^{i\phi_{-}(t,x,\xi)}a_-(t,x,\xi)\hat{f}(\xi)\,d\xi\delta(|x-R\boldsymbol{\theta}|^2-r^2)\,dx.
\end{gather*}  
We drop subscripts in the integral for now and consider only $M_+$,

$$
M_+f(t,\theta) = \frac{1}{(2\pi)^3}\iint_{\R^2\times\R^2} e^{i\phi (t,x,\xi)}a(t,x,\xi) \hat{f}(\xi)\,d\xi \delta(\left| x - R\boldsymbol{\theta}\right|^2 - r^2)\,dx.$$
We make use of the fact that $\frac{1}{2\pi} \int e^{i\lambda(\left| x - R\boldsymbol{\theta}\right|^2 - r^2)}d\lambda = \delta(\left| x - R\boldsymbol{\theta}\right|^2 - r^2)$ to say

$$
M_+f(t,\theta) = \frac{1}{(2\pi)^4}\iiint_{\R^2\times\R\times\R^2} e^{i\phi (t,x,\xi) + i\lambda(\left| x - R\boldsymbol{\theta}\right|^2 - r^2) }a(t,x,\xi) \hat{f}(\xi)\,d\xi\,d\lambda \,dx.
$$
Lastly, we unpack the Fourier transform of $f$ to get

$$
M_+f(t,\theta) = \frac{1}{(2\pi)^4}\iiiint_{\R^7} e^{i\phi (t,x,\xi) + i\lambda(\left| x - R\boldsymbol{\theta}\right|^2 - r^2) - iy\cdot\xi }a(t,x,\xi) f(y)dy\,d\xi\,d\lambda\,dx,$$
where we have identified $\R^2\times\R^2\times\R\times\R^2$ with $\R^7$.  This is an indication that the measurement operator $M_+$ is a Fourier Integral Operator with phase function
$$
\Phi(t,\theta,y;\lambda,x,\xi) = \phi (t,x,\xi) + \lambda(\left| x - R\boldsymbol{\theta}\right|^2 - r^2) - y\cdot\xi.
$$
It can be shown that this phase function is non-degenerate in the sense of \cite{taylor81}.

One issue with this phase function is that $\Phi$ is not homogeneous of degree one in the fiber variables $(\lambda,x,\xi)$, but this can be fixed by making a change of variable.  Let $\tilde{x} := x \left| (\xi, \lambda )\right|$  and define $\tilde{\Phi}(t,\theta,y;\lambda, \tilde{x},\xi) := \Phi(t,\theta,y;\lambda,\frac{\tilde{x}}{\left| (\xi,\lambda )\right|},\xi)$.  This makes $\tilde{\Phi}$ homogeneous of degree one in the variables $(\lambda, \tilde{x}, \xi)$, and so we can proceed.  After making this change of variable, we now write $M_+$ as $$M_+f = \frac{1}{(2\pi)^4}\iiiint_{\R^7} e^{\tilde{\Phi}(t,\theta,y;\lambda,\tilde{x},\xi)}\tilde{a}(t,\tilde{x},\xi,\lambda)f(y)\,dy\,d\xi\,d\lambda\,d\tilde{x}$$
where $$\tilde{a}(t,\tilde{x},\xi,\lambda) = a\left(t,\frac{\tilde{x}}{|(\xi,\lambda)|},\xi\right) |(\xi,\lambda)|^{-2}.$$  Also, because $a$ is an amplitude of order $0$, by \cite{hormander1971} we know that $\tilde{a}(t,\tilde{x},\lambda,\xi)$ is an amplitude of order $-1/2$.  Note that this change of variable does not af{}fect the characteristic manifold for $\Phi$, for
\begin{align*}
\tilde{\Phi}_\lambda &= \Phi_\lambda + \Phi_x \left( \frac{\tilde{x}}{|(\xi,\lambda)|}\right)_\lambda,\\
\tilde{\Phi}_{\tilde{x}} &= \Phi_x \frac{1}{|(\xi,\lambda)|}, \\
\tilde{\Phi}_\xi &= \Phi_\xi + \Phi_x \left( \frac{\tilde{x}}{|(\xi,\lambda)|}\right)_\xi,
\end{align*}
so that $\tilde{\Phi}_{\tilde{x}} = 0$ if and only if $\Phi_{x} = 0$, from which it is clear that $\tilde{\Phi}_{\lambda,\tilde{x},\xi} = 0$ if and only if $\Phi_{\lambda,x,\xi} = 0$.

The characteristic manifold of this FIO is defined as the set $\Sigma_+ = \left\{ (t,\theta,y;\lambda,x,\xi) \mid \Phi_{\lambda,x,\xi} = 0\right\}$.  Taking the derivative of $\Phi$ with respect to $(\lambda, x,\xi)$, we see that this gives the system of equations

\begin{align*}
\left|x - R\boldsymbol{\theta}\right|^2 - r^2 &= 0 \\
\phi_x + 2\lambda  (x - R\boldsymbol{\theta}) &= 0 \\
\phi_\xi - y &= 0 \\
\end{align*}

We've made no assumption on which experimental setup we've chosen to examine so far.  There are two dif{}ferent cases, as mentioned above: (1)$R-r \ge 1$, which we will call the small radius case, and (2)$R=1$ with $r\ge 2$, which we will call the large radius case. The analysis of these two cases are largely the same, but with a few key dif{}ferences.  We examine both cases.
\subsection{Case 1: Small Radius}
From the system of equations obtained by looking at the characteristic manifold of the FIO, we see that $\phi_\xi = y$, and by the geometric optics construction, $x$ lies on the geodesic $\gamma_{y,\hat{\xi}}(t)$ issued from $(y,\hat{\xi})$ where $\hat{\xi} = \xi/(c|\xi|)$ is the unit covector in the metric identified with a unit vector, and $(\gamma_{y,\hat{\xi}}(t),c|\xi|\dot{\gamma}_{y,\hat{\xi}}(t)) = (x,\phi_x)$.  Now from the first equation, we know that $x = (x_1,x_2)$ must lie on the circular integrating detector  of radius $r$ with center $R\boldsymbol{\theta}$. So $x$ is the intersection of the geodesic $\gamma_{y,\hat{\xi}}$ with the circle defined by $|x - R\boldsymbol{\theta}|^2 = r^2$.  There are two of these points of intersection in general (we will show that the geodesic does not intersect the circular integrating detector tangentially if the singularity is to be detected), which we label $x_{+,1}(y,\xi)$ and $x_{+,2}(y,\xi)$.  Also denote the times at which these intersections occur $t_{+,1}(y,\xi)>0$ and $t_{+,2}(y,\xi)>0$ respectively.

Finally, we have from the second equation $\phi_x = 2\lambda (R\boldsymbol{\theta} - x)$.  This tells us that $\dot{\gamma}_{y,\xi}(t)$ is parallel to $x - R\boldsymbol{\theta}$, provided that $\lambda \neq 0$ (in which case $\gamma_{y,\hat{\xi}}(t)$ intersects the circular integrating detector tangentially).  Supposing for a moment that $\lambda = 0$, then taking magnitudes on both sides of the equation $\phi_x(t,x,\xi) = 2\lambda(R\boldsymbol{\theta} - x)$, we obtain $|\phi_x(t,x,\xi)| = 0\implies c(y)|\xi||\dot{\gamma}_{y,\hat{\xi}}(t)| = 0$.  This in turn means $|\dot{\gamma}_{y,\hat{\xi}}(t)| = 0$ as $c(y) \neq 0$ for any $y\in \R^2$ and $\xi$ is a non zero vector.  But we know near the integrating detectors, that $|\dot{\gamma}_{y,\hat{\xi}}(t)| = 1$, a contradiction, so $\lambda \neq 0$ on the characteristic manifold. And so the geodesic intersects the circular detector perpendicularly.  We know that $\phi_x = c|\xi|\dot{\gamma}_{y,\hat{\xi}}(t)$, And $c\equiv 1$ outside of the region of interest, so that $|\lambda_{+,i}| = c(y)/(2r)|\xi||(\dot{\gamma}_{y,\hat{\xi}}(t_{+,i}))| = c(y)/(2r)|\xi|$.  Note that $\lambda_{+,1} > 0$ and $\lambda_{+,2} < 0$, so that $\lambda_{+,1} = c(y)/(2r)|\xi|$ and $\lambda_{+,2} = -\lambda_{+,1}$.  We'll simply denote $\lambda_+ := \lambda_{+,1}$.  Because $\gamma_{y,\hat{\xi}}(t)$ intersects the circular integrating detector perpendicularly, it must go through the center of the circular integrating detector, as outside of $\Omega$, we know that $c \equiv 1$ implies that $\gamma_{y,\hat{\xi}}$ is a straight line near the integrating detectors, and so we see that $\boldsymbol{\theta} = \gamma_{y,\hat{\xi}}(t_{+,1}(y,\xi) + r)/R = \gamma_{y,\hat{\xi}}(t_{+,2}(y,\xi) - r)/R$, where $r$ is the fixed radius of the circular integrating detector.  This gives then the entire characteristic manifold parametrized by $(y,\xi)$, giving a smooth manifold of dimension 4 consisting of 2 connected parts.  Define $\Sigma_{+,1} = \left\{ (t_{+,1}(y,\xi),\boldsymbol{\theta}(y,\xi),y;\lambda_+(y,\xi),x_{+,1}(y,\xi),\xi) \mid (y,\xi) \in T^*(\Omega) \setminus \left\{ 0 \right\}\right\}$ and $\Sigma_{+,2} = \left\{ (t_{+,2}(y,\xi),\boldsymbol{\theta}(y,\xi),y;-\lambda_+(y,\xi),x_{+,2}(y,\xi),\xi) \mid (y,\xi) \in T^*(\Omega) \setminus \left\{ 0 \right\}\right\}$.  Then $\Sigma_+ = \Sigma_{+,1} \cup \Sigma_{+,2}$ as a disjoint union.

We now write down the canonical relation given by
$$
\Sigma_+ \ni (t,\boldsymbol{\theta},y;\lambda,x,\xi) \mapsto (t,\boldsymbol{\theta},y;\Phi_t,\Phi_{\theta},\Phi_y)
$$
This mapping can be calculated as

$$
\Sigma_+ \ni (t,\boldsymbol{\theta},y;\lambda,x,\xi) \mapsto (t,\boldsymbol{\theta},y;-c(y)|\xi|,-\frac{R}{r}c(y)|\xi|(x-(x\cdot \boldsymbol{\theta})\boldsymbol{\theta}),-\xi)
$$

The analysis for $\sigma = -$ is the same giving us $t_{-,i}(y,\xi) < 0,x_{-,i}(y,\xi)$ for $i=1,2$ and $\boldsymbol{\theta}_-$.  We also see that $\lambda_{-,1} < 0$ and $\lambda_{-,2} > 0$.  Denote $\lambda_{-}:= \lambda_{-,1}$. We see that the maps

$$
\Sigma_{\pm} \ni (t,\boldsymbol{\theta},y;\lambda,x,\xi)\mapsto (t,\boldsymbol{\theta},y;\mp c(y)|\xi|,\mp \frac{R}{r}c(y)|\xi|(x-(x\cdot\boldsymbol{\theta})\boldsymbol{\theta}),-\xi)
$$
are smooth and of full rank, and so the canonical relations associated to the operators $M_+$ and $M_-$ are given by

$$
C_{\pm} := \left\{(t,\boldsymbol{\theta},\Phi_t,\Phi_{\theta};y,\xi) \mid (t,\boldsymbol{\theta},y;\lambda,x,\xi) \in \Sigma_\pm \right\}$$
$$
=\left\{ \left(t_{\pm,i}(y,\xi), \boldsymbol{\theta}_{\pm}(y,\xi), \mp c(y)|\xi|, \mp \frac{R}{r}c(y)|\xi|(\gamma_{y,\hat{\xi}}(t_{\pm,i}(y,\xi))-(x_{\pm,i}\cdot\boldsymbol{\theta}_{\pm})\boldsymbol{\theta}_\pm(y,\xi));y,\xi\right)\right\}$$
such that $(y,\xi)\in T^*(\Omega)\setminus{0}$, where $i = 1,2$.  Here, $\boldsymbol{\theta}_+ = \gamma_{y,\hat{\xi}}(t_{+,1}(y,\xi) + r)/R$ and $\boldsymbol{\theta}_- = \gamma_{y,\hat{\xi}}(t_{-,1}(y,\xi) - r)/R$.  Writing this as a mapping, we have

\begin{gather*}
{C_{\pm} : (y,\xi)} \mapsto 
\left\{ {\begin{matrix}
\left(t_{\pm,1}, \frac{1}{R}\gamma_{y,\hat{\xi}}(t_{\pm,1} \pm r), \mp c(y)|\xi|, \pm \frac{R}{r}c(y)|\xi|\left(\left(x_{\pm,1}\cdot\theta_\pm\right)\gamma_{y,\hat{\xi}}(t_{\pm,1} \pm r)-\gamma_{y,\xi}(t_{\pm,1})\right)\right) \\ \\
\left(t_{\pm,2}, \frac{1}{R}\gamma_{y,\hat{\xi}}(t_{\pm,2} \mp r), \mp c(y)|\xi|, \mp \frac{R}{r}c(y)|\xi|\left(\left(x_{\pm,2}\cdot\theta_{\pm}\right)\gamma_{y,\hat{\xi}}(t_{\pm,2}\mp r)-\gamma_{y,\xi}(t_{\pm,2})\right)\right)
\end{matrix}}\right.
\end{gather*}
Where $x_{\pm,i} = x_{\pm,i}(y,\xi)$, $\theta_{\pm} = \theta_{\pm}(y,\xi)$, and $t_{\pm,1} = t_{\pm,1}(y,\xi)$ and similarly, $t_{\pm,2} = t_{\pm,2}(y,\xi)$. The canonical relations for the operators $M_+$ and $M_-$ are each one to two and of the above form.  The above has shown the following:

\begin{theorem}\label{thm:small}
For $R-r \ge 1$, the operator $Mf = M_+f + M_-f$ defined above is a Fourier Integral Operator of order $-1/2$ with canonical relation given by
\begin{gather*}
    C = C_+ \cup C_-
\end{gather*}
where
\begin{gather*}
C_{\pm} =\left\{ \left(t_{\pm,i}, \boldsymbol{\theta}_{\pm}, \mp c(y)|\xi|, \pm (-1)^{i}\frac{R}{r}c(y)|\xi|\left(\gamma_{y,\hat{\xi}}(t_{\pm,i})-(x_{\pm,i}(y,\xi)\cdot\boldsymbol{\theta}_{\pm})\boldsymbol{\theta}_{\pm}\right);y,\xi\right)\right\},
\end{gather*}
($t_{\pm,i} = t_{\pm,i}(y,\xi)$ and $\boldsymbol{\theta}_{\pm} = \boldsymbol{\theta}_{\pm}(y,\xi)$) with $(y,\xi)\in T^*(\Omega)\setminus{0}$ and $i=1,2$.
\end{theorem}
Note that this canonical relation is locally one to four.
\subsection{Case 2: Large Radius}
In this case, the analysis is almost entirely the same, except there is only one point of intersection of the geodesic $\gamma_{y,\hat{\xi}}(t)$ with the circular integrating detector defined by $|x-\boldsymbol{\theta}| = r$ (see Figure \ref{fig:cases}).  We then have $\Sigma_{\pm} = \{(t_{\pm}(y,\xi),\boldsymbol{\theta}_{\pm}(y,\xi),y;\pm\lambda(y,\xi),x_{\pm}(y,\xi),\xi)\mid (y,\xi) \in T^*(\Omega)\setminus\{0\}\}$, and the canonical relations are given by
$$
C_{\pm} =\left\{ ((t_{\pm}, \boldsymbol{\theta}_{\pm}(y,\xi), \mp c(y)|\xi|, \mp \frac{1}{r}c(y)|\xi|(\gamma_{y,\hat{\xi}}(t_{\pm})-(x_{\pm}\cdot\boldsymbol{\theta}_{\pm})\boldsymbol{\theta}_\pm(y,\xi));y,\xi) \mid (y,\xi)\in T^*(\Omega)\setminus{0}\right\}.$$
Here $t_\pm = t_\pm(y,\xi)$. The canonical relations are each locally one to one in this case, and we have an analogous result as in the first case:
\begin{theorem}\label{thm:large}
For $R = 1$ and $r\ge 2$, the operator $Mf = M_+f + M_-f$ defined above is a Fourier Integral Operator of order $-1/2$ with canonical relation given by
\begin{gather*}
    C = C_+ \cup C_-
\end{gather*}
where
\begin{gather*}
C_{\pm} = \left\{ ((t_{\pm}, \boldsymbol{\theta}_{\pm}(y,\xi), \mp c(y)|\xi|, \mp \frac{1}{r}c(y)|\xi|(\gamma_{y,\hat{\xi}}(t_{\pm})-(x_{\pm}\cdot\boldsymbol{\theta}_{\pm})\boldsymbol{\theta}_\pm(y,\xi));y,\xi) \mid (y,\xi)\in T^*(\Omega)\setminus{0}\right\},
\end{gather*}
with $t_\pm = t_\pm(y,\xi)$.
\end{theorem}
This canonical relation is locally one to two, as each individual canonical map is locally one to one.

Note that the maps $C_+$ and $C_-$ are not globally one to 1, although each are locally one to one, for suppose (looking only at $C_+$ for a moment) $(t(y_1,\xi_1),\boldsymbol{\theta}(y_1,\xi_1), \tau(y_1,\xi_1), \omega(y_1,\xi_1)) = (t(y_2,\xi_2),\boldsymbol{\theta}(y_2,\xi_2), \tau(y_2,\xi_2), \omega(y_2,\xi_2))$ where $\tau(y,\xi) = - c(y)|\xi|$ and $\omega(y,\xi) = -\frac{1}{r}c(y)|\xi|(x -(x \cdot\boldsymbol{\theta})\boldsymbol{\theta}(y,\xi))$ and $x = \gamma_{y,\hat{\xi}}(t(y,\xi))$.  We'll call $x_i = x(y_i,\xi_i)$ for $i=1,2$.  Then clearly we have $t(y_1,\xi_1) = t(y_2,\xi_2)$ and $\boldsymbol{\theta}(y_1,\xi_1) = \boldsymbol{\theta}(y_2,\xi_2)$, which we'll just label $t$ and $\boldsymbol{\theta}$ respectively.  We also clearly have $c(y_1)|\xi_1| = c(y_2)|\xi_2|$.  Suppose for a moment that $(x_1 - x_2)\cdot \boldsymbol{\theta}^{\bot}=0$ where $\boldsymbol{\theta}^\bot = (-\sin(\theta),\cos(\theta))$ is the unit vector perpendicular to $\boldsymbol{\theta}$.  Then note that
\begin{gather*}
    \omega(y_1,\xi_1) - \omega(y_2,\xi_2) = x_1 - (x_1\cdot \boldsymbol{\theta})\boldsymbol{\theta} - (x_2 - (x_2\cdot \boldsymbol{\theta})\boldsymbol{\theta}) \\
    = (x_1 - x_2) - [(x_1 - x_2)\cdot \boldsymbol{\theta}]\boldsymbol{\theta}
\end{gather*}
Now note that $(\omega(y_1,\xi_1) - \omega(y_2,\xi_2))\cdot \boldsymbol{\theta} = 0$ and $(\omega(y_1,\xi_1) - \omega(y_2,\xi_2))\cdot \boldsymbol{\theta}^\bot = 0$. Because $\boldsymbol{\theta}$ and $\boldsymbol{\theta}^\bot$ are linearly independent, this shows that $\omega(y_1,\xi_1) - \omega(y_2,\xi_2) = 0$ or in other words that $\omega(y_1,\xi_1) = \omega(y_2,\xi_2)$.  This shows that, provided $t(y_1,\xi_1) = t(y_2,\xi_2)$, $\boldsymbol{\theta}(y_1,\xi_1) = \boldsymbol{\theta}(y_2,\xi_2)$, $c(y_1)|\xi_1| = c(y_2)|\xi_2|$ and $(x_1 - x_2)\cdot \boldsymbol{\theta}^\bot = 0$, that $(y_1,\xi_1)$ and $(y_2,\xi_2)$ get mapped to the same point under the canonical relation. However, for $c$ close enough to $1$, this will not happen locally.

\section{Injectivity}

\subsection{Case 1: Small Radius}
Let $u(t,x)$ be the solution to (\ref{eq:1}) and $\Gamma \subset S^1$ be open.  In local coordinates, suppose $\Gamma$ is the open interval given by $\Gamma = (\theta_1, \theta_2)$, with $0 \le \theta_1 < \theta_2 < 2\pi$.  Then, for $\boldsymbol{\theta} = (\cos(\theta),\sin(\theta)) \in \Gamma$, we have 
$$Mf(t,\theta) = \frac{1}{2\pi} \int_0^{2\pi}u(t,R_0\cos(\theta) + r\cos(\alpha),R_0\sin(\theta) + r\sin(\alpha))\,d\alpha$$
Where $R_0$ is the distance from the orgin to the center of the circular integrating detector. For fixed $r>0$, we may view $R$ as variable with $R \ge R_0 \ge 1+r$ (i.e. we may translate the circular integrating detectors away from the region of interest).  We denote this by letting $Mf(t,\theta)$ vary with $R$ and denote the operator then as $Mf(t,\theta;R)$ and note $Mf(t,\theta) = Mf(t,\theta;R_0)$.  Let $P(t,\theta,R) = Mf(t,\theta;R)$.  Then we see that
\begin{equation}\label{eq:2}
\begin{cases}R^2 P_{tt}(t,\theta,R) - R(R P_R(t,\theta,R))_R - P_{\theta\theta}(t,\theta,R) = 0 & (t,\theta,R) \in \mathbb{R}^+\times \Gamma \times [R_0,\infty),\\
 P(t,\theta,R_0) = Mf(t,\theta)\mid_{\theta \in \Gamma}, & (t,\theta) \in \R^+\times \Gamma, \\
 P(0,\theta,R) = 0, & (\theta,R) \in \Gamma \times [R_0,\infty),\\
 P_t(0,\theta,R) = 0, & (\theta,R)\in \Gamma \times [R_0,\infty).
\end{cases}
\end{equation}
This can be seen as follows:  first note to save space, $u$ and all of its partial derivatives are understood to be evaluated at $(t,R\cos(\theta) + r\cos(\alpha),R\sin(\theta) + r\sin(\alpha))$.  We have

$$P_\theta(t,\theta,R) = \frac{1}{2\pi}\int_0^{2\pi} (-R\sin(\theta))u_x  + (R\cos(\theta))u_y \,d\alpha,$$
and so we have $P_{\theta \theta}$ can be written:
$$
\frac{1}{2\pi}\int_0^{2\pi}R^2\sin^2(\theta)u_{xx} - 2R^2\sin(\theta)\cos(\theta)u_{xy} + R^2\cos^2(\theta)u_{yy}-R\cos(\theta)u_x - R\sin(\theta)u_y\,d\alpha.$$
Noting that 
\begin{align*}
 R^2\sin^2(\theta)u_{xx} + R^2\cos^2(\theta)u_{yy} &= R^2\sin^2(\theta)u_{xx} + R^2\cos^2(\theta)u_{xx} + R^2\sin^2(\theta)u_{yy}+R^2\cos^2(\theta)u_{yy} \\
 &- (R^2\cos^2(\theta)u_{xx} + R^2\sin^2(\theta)u_{yy})\\
 &= R^2\Delta u - (R^2\cos^2(\theta)u_{xx} + R^2\sin^2(\theta)u_{yy}),
\end{align*}
we see that $P_{\theta \theta}(t,\theta,R)$ is given by

$$
\frac{1}{2\pi}\int_0^{2\pi}R^2\Delta u - (R^2\cos^2(\theta)u_{xx} + R^2\sin^2(\theta)u_{yy} + 2R^2\sin(\theta)\cos(\theta)u_{xy} + R\cos(\theta)u_x + R\sin(\theta)u_y)\,d\alpha.$$

Remembering that $\Delta u$ above is evaluated at $(t,R\cos(\theta)+r\cos(\alpha),R\sin(\theta)+r\sin(\alpha)),$ where $c \equiv 1$, we have $P_{\theta \theta}(t,\theta,R)=$

$$
\frac{1}{2\pi}\int_0^{2\pi}R^2 c^2 \Delta u - (R^2\cos^2(\theta)u_{xx} + R^2\sin^2(\theta)u_{yy} + 2R^2\sin(\theta)\cos(\theta)u_{xy} + R\cos(\theta)u_x + R\sin(\theta)u_y)\,d\alpha,$$
so that

$$
P_{\theta \theta}(t,\theta,R) = \frac{1}{2\pi}R^2 \int_0^{2\pi} u_{tt}\,d\alpha - \frac{1}{2\pi}\int_0^{2\pi}(R^2\cos^2(\theta)u_{xx} + R^2\sin^2(\theta)u_{yy} + 2R^2\sin(\theta)\cos(\theta)u_{xy} $$$$ + R\cos(\theta)u_x + R\sin(\theta)u_y)\,d\alpha,
$$
and $$\frac{1}{2\pi}\int_0^{2\pi}u_{tt}\,d\alpha = \frac{1}{2\pi}\left(\int_0^{2\pi}u\,d\alpha\right)_{tt} = P_{tt}(t,\theta,R).$$
So we only need to show that $R(RP_R)_R=$

$$
\frac{1}{2\pi} \int_0^{2\pi}(R^2\cos^2(\theta)u_{xx} + R^2\sin^2(\theta)u_{yy} + 2R^2\sin(\theta)\cos(\theta)u_{xy} + R\cos(\theta)u_x + R\sin(\theta)u_y)\,d\alpha.
$$
This follows from direct calculation and the chain rule.

Rewriting (\ref{eq:2}) in a more standard form, we have:
\begin{equation}\label{eq:3}
\begin{cases} P_{tt}(t,\theta,R) - \frac{1}{R}(R P(t,\theta,R)_R)_R - \frac{1}{R^2} P_{\theta\theta}(t,\theta,R) = 0, & (t,\theta,R) \in \mathbb{R}^+\times \Gamma \times [R_0,\infty),\\
 P(t,\theta,R_0) = Mf(t,\theta)\mid_{\theta \in \Gamma}, & (t,\theta) \in \R^+ \times \Gamma, \\
 P(0,\theta,R) = 0, & (\theta,R) \in \Gamma \times [R_0,\infty)\\
 P_t(0,\theta,R) = 0, & (\theta,R)\in \Gamma \times [R_0,\infty)
\end{cases}
\end{equation}

(\ref{eq:3}) simply says that $P(t,\theta,R)$ is a solution to the constant speed wave equation in polar coordinates with initial condition $P(0,\theta,R) = 0$ for $R \ge R_0 > 1+r$ and $P(t,\theta,R_0) = Mf(t,\theta)$ for $t\ge 0$ and $\theta \in \Gamma$.  It is well known that the wave equation on a Riemannian manifold has a unique solution, and so the solution to (\ref{eq:3}) is unique.  To show then that $Mf(t,\theta)\mid_{\theta \in \Gamma}$ is uniquely determined by $f$, we need to show that $Mf(t,\theta)\mid_{\theta \in \Gamma} = 0 \implies f(x) = 0$ by the linearity of $M$.  So then, we assume that $Mf(t,\theta) \equiv 0$ for $t \in \R^+$ and $\theta \in \Gamma$.  By the uniqueness of solutions to the wave equation, this tells us that $P(t,\theta, R) \equiv 0$ for any $(t,\theta, R) \in \R^+ \times \Gamma \times [R_0 ,\infty)$.  Let $T \in \R$.  We may extend $u(t,x)$ in an even way for $|t| < T$ so that $u$ is still a solution to the wave equation, and so we may extend $P(t,\theta,R)$ in an even way such that $P(t,\theta,R) = 0$ for $|t| < T$. By finite speed of propagation, we know that $\textnormal{supp }u(T,\cdot) \subset B_{1+T}(0)$.  Let $\theta_0 \in \Gamma$.  Note that the set $\mathcal{A} = \{ (\theta, R)\mid \theta \in \Gamma,\,R > R_0\}$ is open and connected in $\R^2$, because $\Gamma \subset S^1$ is an open interval. We know that the circular Radon transform of $u(T,\cdot)$ is 0 for any $\theta \in \Gamma$ and for any $R > R_0$.  Further, because the interior of these circular integrating detectors is just $D(\theta,R) = B_r(R\boldsymbol{\theta})$, we can take $R > R_0$ large enough so that $\textnormal{supp }u(T,\cdot) \cap D(\theta_0,R) = \emptyset$, because $\textnormal{supp }u(T,\cdot)$ is contained in a bounded set.  It follows then by Theorem 1.2 in \cite{quinto94} that $\textnormal{supp }u(T,\cdot)$ is disjoint from $\cup_{(\theta,R) \in \mathcal{A}} D(\theta,R)$.  So in particular, there is a neighborhood $V$ of $x_0 = R_0\boldsymbol{\theta}_0$ such that $u(T,x) = 0$ on $V$.  $T$ was chosen arbitrarily, so this result holds for all $|t| < T$, and so by Tataru's unique continuation, $u(t,x) \equiv 0$ in the domain of influence $|t| + \textnormal{dist}(x,x_0) < T$.  So, taking $T$ large enough so that $\textnormal{dist}(x,x_0) < T$ for all $x \in \Omega$, we see that $u(0,x) = f(x) \equiv 0$, and so $Mf(t,\theta)\mid_{[0,T] \times \Gamma}$ is uniquely determined by $f$.

\subsection{Case 2: Large Radius}
Again, we consider $u(t,x)$ a solution to (\ref{eq:1}).  We consider only the full data case $\Gamma = S^1$. Then in this case, in which $R = 1$, we have that the measurement operator $Mf(t,\theta)$ is given by
$$
Mf(t,\theta) = \frac{1}{2\pi}\int_0^{2\pi} u(t,\cos(\theta) + r_0\cos(\alpha),\sin(\theta) + r_0\sin(\alpha))\,d\alpha
$$
where $r_0 \ge 2$ is the fixed radius of the circular integrating detector.  We may however view $r$ as variable, noting the operator with variable detector radius by 

$$Mf(t,\theta;r) = \frac{1}{2\pi}\int_0^{2\pi} u(t,\cos(\theta) + r\cos(\alpha),\sin(\theta) + r\sin(\alpha))\,d\alpha$$
Where $r \ge r_0 \ge 2$.
Let $P(t,\theta,r) = Mf(t,\theta;r)$, where $\theta \in \Gamma$ and $r \ge r_0 > 2$.  Note that $P(t,\theta,r_0) = Mf(t,\theta)$ for $\theta \in \Gamma$.  It follows then that $P(t,\theta,r)$ solves the following PDE:
\begin{equation}\label{eq:4}
 \begin{cases}
  P_{tt}(t,\theta,r) - \frac{1}{r}(rP_{r}(t,\theta,r))_r = 0 & (t,\theta,r) \in \mathbb{R}^+ \times \Gamma \times [r_0,\infty) \\
  P(t,\theta,r_0) = Mf(t,\theta) & (t,\theta) \in \mathbb{R}^+ \times \Gamma \\
  P(0,\theta,r) = 0 & (\theta,r)\in \Gamma \times [r_0,\infty) \\
  P_t(0,\theta,r) = 0 & (\theta,r) \in \Gamma \times [r_0,\infty)
 \end{cases}
\end{equation}
We see this as follows:  (Again note that $u$ and all of its partial derivatives are understood to be evaluated at $(t,\cos(\theta) + r\cos(\alpha),\sin(\theta) + r\sin(\alpha))$.)  We have
$$
P_r(t,\theta,r) = \frac{1}{2\pi}\int_0^{2\pi} \cos(\alpha)u_x + \sin(\alpha)u_y\,d\alpha.
$$
We integrate by parts to get

$$
P_r(t,\theta,r) = \frac{1}{2\pi}\int_0^{2\pi} (r\sin^2 \alpha) u_{xx} -2r\cos(\alpha)\sin(\alpha)u_{xy} + (r\cos^2(\alpha))u_{yy}\,d\alpha.
$$
Then, we use the fact that
$$
r\sin^2(\alpha)u_{xx} + r\cos^2(\alpha)u_{yy} = r \Delta u - r\cos^2(\alpha) u_{xx} - r\sin^2(\alpha)u_{yy},
$$
to obtain

$$
P_r(t,\theta,r) = \frac{1}{2\pi}\int_0^{2\pi} r\Delta u - r(\cos^2(\alpha)u_{xx} + 2\cos(\alpha)\sin(\alpha)u_{xy} + \sin^2(\alpha)u_{yy})\,d\alpha.$$
We recall that $u$ and its derivatives are evaluated at $(t,\cos(\theta) + r\cos(\alpha),\sin(\theta) + r\sin(\alpha))$ where $c \equiv 1$, so that $\Delta u = u_{tt}$ there and we see that

\begin{align*}
P_r(t,\theta,r) &= r \frac{1}{2\pi}\int_0^{2\pi} u_{tt}\,d\alpha - r\frac{1}{2\pi}\int_0^{2\pi}\cos^2(\alpha)u_{xx} + 2\cos(\alpha)\sin(\alpha)u_{xy} + \sin^2(\alpha)u_{yy}\,d\alpha \\
&= rP_{tt}(t,\theta,r) - r\frac{1}{2\pi}\int_0^{2\pi}\cos^2(\alpha)u_{xx} + 2\cos(\alpha)\sin(\alpha)u_{xy} + \sin^2(\alpha)u_{yy}\,d\alpha \\
&= rP_{tt}(t,\theta,r) - rP_{rr}(t,\theta,r).
\end{align*}
The last line can be seen by direct calculation of $P_{rr}(t,\theta,r)$.  Rearranging we see $rP_{tt} = rP_{rr} + P_r$ and (\ref{eq:4}) then follows.  This partial dif{}ferential equation is the wave equation with axial symmetry, and so solutions to this equation are again unique.  

To show that $Mf(t,\theta)$ determines $f$ uniquely for $\theta \in \Gamma$, by the linearity of $M$, we need only show that $Mf(t,\theta)= 0$ for all $(t,\theta) = [0,T)\times \Gamma \implies f = 0$.  We note that $Mf(t,\theta) = 0$ for all $(t,\theta) \in [0,T) \times \Gamma \implies P(t,\theta,r) \equiv 0$, because $P(t,\theta,r) \equiv 0$ clearly solves the PDE (\ref{eq:4}) in this case, and this solution is unique.  For any $T\in \R$, we know that $\textnormal{supp } u(T,\cdot) \subset B_{1+T}(0)$, and so for any finite $T$, $u(T,\cdot)$ has bounded support.  Let $\theta_0 \in \Gamma$, $T\in \R^+$.  Let $r_1 = \sup \{ r \ge r_0 \mid C(\theta_0,r)\cap \textnormal{supp } u(T,\cdot) \ne \emptyset\}$, where $C(\theta_0, r)$ is the circle centered at $\boldsymbol{\theta}_0$ with radius $r$. We know, because $\supp u(T,\cdot) \subset B_{1+T}(0)$, that $\{ r\ge r_0 \mid C(\theta_0,r) \cap \supp u(T,\cdot) \ne \emptyset\}$ is bounded above by $2 + T$.  Assume that this set is nonempty so that $r_1 \ge r_0$ is finite.  Then we have $C(\theta_0,r_1) \cap \supp u(T,\cdot) \ne \emptyset$ by the compactness of $\supp u(T,\cdot)$.  Let $x \in C(\theta_0,r_1)\cap \supp u(T,\cdot)$, and let $(x,\xi) \in N^*(C(\theta_0,r_1))$.  By construction, $\supp u(T,\cdot)$ is on one side of $C(\theta_0,r_1)$ at $x$, so that by Theorem 8.5.6 in \cite{hormander2009}, we have that $(x,\xi) \in \textnormal{WF}_A(u(T,\cdot))$.  Note that as in the terminology of \cite{quinto94}, it is impossible for $x$ to be a $C(\theta_0,r_1)$ self mirror point, for tracing $(x,\xi)$ back along the geodesic defined by $(x,\xi)$, we see that the geodesic never intersects the interior of $B_1(0)$, which is impossible.  So then, there are two cases we must consider.  First we consider the case where the $C(\theta_0,r_1)$ mirror point of $x$, which we will call $\tilde{x}$, is not in the intersection $C(\theta_0,r_1) \cap \supp u(T,\cdot)$.  Then, by the compactness of $\supp u(T,\cdot)$, we have $u(T,\cdot) = 0$ in a neighborhood of $\tilde{x}$.  We also have $P(T,\theta,r) = 0$ in a neighborhood of $(\theta_0,r_1) \in \Gamma \times [r_0,\infty)$.  It then follows  by Proposition 2.4 of \cite{quinto94} that $(x,\xi) \not\in \textnormal{WF}_A(u(T,\cdot))$, a contradiction to $(x,\xi)\in \textnormal{WF}_A(u(T,\cdot))$.  It follows that $\{r \ge r_0 \mid C(\theta_0,r) \cap \supp u(T,\cdot) \ne \emptyset \}$ is empty and so $\supp u(T,\cdot) \subset \bar{B}_{r_0}(\boldsymbol{\theta}_0)$.

The second case we consider is that $\tilde{x}$ is in the intersection $C(\theta_0,r_1)\cap \supp u(T,\cdot)$.  We'll show that there then exists $\theta_1 \in S^1$ and $r_2 > r_1$, such that $C(\theta_1,r_2) \cap \supp u(T,\cdot) = \{ \tilde{x} \}$.  Assume for now that this is the case, and let $(\tilde{x},\tilde{\xi}) \in N^*(C(\theta_1,r_2))$.  Then again, we have by construction that $\supp u(T,\cdot)$ is on one side of $C(\theta_1,r_2)$ at $\tilde{x}$ and so $(\tilde{x},\tilde{\xi}) \in \textnormal{WF}_A(u(T,\cdot))$.  It also follows as before that $u(T,\cdot)$ is zero in a neighborhood of the $C(\theta_1,r_2)$ mirror point of $\tilde{x}$.  We then see again from Proposition 2.4 of \cite{quinto94} that $(\tilde{x},\tilde{\xi})\not\in \textnormal{WF}_A(u(T,\cdot))$, a contradiction.  It follows again that $\supp u(T,\cdot) \subset \bar{B}_{r_0}(\boldsymbol{\theta}_0)$.
That $f$ is zero then follows from Tataru's unique continuation as in the small radius case.

Now we show the existence of the circle $C(\theta_1,r_2)$ with property that $C(\theta_1,r_2)\cap \supp u(T,\cdot) = \{\tilde{x}\}$ mentioned above.  We let $C(\theta_0,r_1)$, and $x,\tilde{x}$ be as above.  We define the sets $C^{-}(\theta_0,r_1) = \{ y \in C(\theta_0,r_1) \mid (y-\boldsymbol{\theta_0})\cdot \boldsymbol{\theta_0} \le 0 \}$, and $C^{+}(\theta_0,r_1) = \{ y \in C(\theta_0,r_1) \mid (y-\boldsymbol{\theta_0})\cdot \boldsymbol{\theta_0} \ge 0\}$.  It is clear that $B_1(0)$ is contained in the interior of the region bounded by $C^{-}(\theta_0,r_1) \cup L(\theta_0,r_1)$, where $L(\theta_0,r_1)$ is the diameter of $C(\theta_0,r_1)$ defined by the vector $\boldsymbol{\theta}_0^{\bot}$.  We may assume without loss of generality that $x \in C^{-}(\theta_0,r_1)$, for if not, we may simply swap the roles of $x$ and $\tilde{x}$ in what follows.  Now $x \in C^{-}(\theta_0,r_1) \implies \tilde{x} \in C^{+}(\theta_0,r_1)$, because $x$ cannot be a $C(\theta_0,r_1)$ mirror point, as we have shown.  The line $\ell(t)$ defined by
\begin{gather*}
    \ell(t) = \boldsymbol{\theta_0} + t(\tilde{x} - \boldsymbol{\theta_0})
\end{gather*}
intersects $S^1$ at 2 points: $\boldsymbol{\theta}_0$ and \begin{gather*}\boldsymbol{\theta}_1 = \boldsymbol{\theta}_0 - \frac{2\boldsymbol{\theta}_0\cdot(\tilde{x} - \boldsymbol{\theta}_0)}{r_1^2}(\tilde{x} - \boldsymbol{\theta}_0).\end{gather*}
This implies there are two distinct circles with centers on $S^1$ such that $(\tilde{x},(\tilde{x} - \boldsymbol{\theta}_0))$ is in the conormal bundle to these circles, namely $C(\theta_0,r_1)$ and $C(\theta_1,r_2)$, where $r_2 = |\tilde{x} - \boldsymbol{\theta}_1|$.  Note that
\begin{gather*}
    r_2 = \left|\tilde{x} - (\boldsymbol{\theta}_0 -\frac{2 (\tilde{x} - \boldsymbol{\theta}_0)\cdot \boldsymbol{\theta}_0}{r_1^2}(\tilde{x} - \boldsymbol{\theta}_0)\right| \\
    = \left| (\tilde{x} - \boldsymbol{\theta}_0) \left( 1 + \frac{2 (\tilde{x}-\boldsymbol{\theta}_0)\cdot \boldsymbol{\theta}_0}{r_1^2}\right) \right| \\
    = r_1 \left| 1 + \frac{2 (\tilde{x} - \boldsymbol{\theta}_0)\cdot\boldsymbol{\theta}_0}{r_1^2}\right| > r_1
\end{gather*}
where the last inequality follows because
\begin{gather*}
    \tilde{x} \in C^+(\theta_0,r_1) \implies \frac{2 (\tilde{x} - \boldsymbol{\theta}_0)\cdot \boldsymbol{\theta}_0}{r_1^2} > 0 \\
    \textnormal{So that }1 + \frac{2 (\tilde{x} - \boldsymbol{\theta}_0)\cdot\boldsymbol{\theta}_0}{r_1^2} > 1,
\end{gather*}
and we see then that $r_2 = r_1 +\frac{2(\tilde{x}-\boldsymbol{\theta}_0)\cdot\boldsymbol{\theta}_0}{r_1}$.
Now, we need only show that $C(\theta_1,r_2) \cap \supp u(T,\cdot) = \{ \tilde{x} \}$.  Clearly by construction, $\tilde{x}$ is in this intersection.  Also note that by the choice of $r_1$, that $C(\theta_0,r) \cap \supp u(T,\cdot) = \emptyset$ for all $r > r_1$.  So let $y \in C(\theta_1,r_2)$.  Then, we have using the triangle inequality
\begin{gather*}
    |y - \boldsymbol{\theta}_0| \ge |y - \boldsymbol{\theta}_1| - |\boldsymbol{\theta}_1 - \boldsymbol{\theta}_0| \\
    = r_2 - \left|\frac{2(\tilde{x}-\boldsymbol{\theta}_0)\cdot\boldsymbol{\theta}_0}{r_1^2}(\tilde{x}-\boldsymbol{\theta}_0)\right| \\
    = r_1,
\end{gather*}
so that $|y - \boldsymbol{\theta}_0| \ge r_1$, and equality holds only when $y - \boldsymbol{\theta}_0 = \alpha(\boldsymbol{\theta}_1 - \boldsymbol{\theta}_0)$ for some $\alpha > 0$, but this is only the case when $y = \tilde{x}$.  So, we've shown that $y \in C(\theta_1,r_2) \implies |y - \boldsymbol{\theta}_0| > r_1$ if $y \neq \tilde{x}$.  In other words, $y \neq \tilde{x} \implies y \in C(\theta_0,r)$ for some $r > r_1$, so that so that $y \not \in \supp u(T,\cdot)$, and in particular, $y \not\in C(\theta_1,r_2)\cap \supp u(T,\cdot)$ if $y \neq \tilde{x}$.  This completes the proof in the large radius case for full data.

\section{Stability}
It is natural to take $f\in H_D(\Omega) \subset H_0^1(\Omega)$ when taking conservation of energy into consideration.
By the above theorems, in both the small and large radius cases, we have that $Mf(t,\theta)$ is an elliptic FIO of order $-1/2$.  We will model finite time measurements in $\left[ 0,T\right]$ by premultiplying $Mf(t,\theta)$ by $\chi \in C_0^\infty ( \mathbb{R} )$ with $\textnormal{supp } \chi \subset [-1,T_1]$ for $T < T_1 < \infty$, and $\chi \equiv 1$ on $[0,T]$.  Now because $M$ is an FIO of order $-1/2$ associated with the graph of the canonical relation $C$, we have $M^*$ is also an FIO of order $-1/2$, and it is associated with the canonical relation $C^{-1}$.  $M^*\chi M$ then is an elliptic pseudodif{}ferential operator of order $-1$.  This implies that a parametrix $B$, which necessarily is an elliptic $\Psi$DO of order 1, exists such that
$$
BM^*\chi M = \textnormal{Id} - R,
$$
where $R$ is a regularizing operator.  We may assume that $B$ is a properly supported $\Psi$DO, which means $B,B^t: \mathcal{E}'(\Omega)\rightarrow \mathcal{E}'(\Omega)$.  From this we see that
$$
\lVert f \rVert_{H^1(\Omega)} \le \lVert BM^* \chi Mf \rVert_{H^1(\Omega)} + \lVert Rf\rVert_{L^2(\Omega)}.
$$
$B$ is a continuous linear operator, so we have
$$
\lVert f \rVert_{H^{1}(\Omega)} \le C \lVert M^*\chi M f \rVert_{H^{2}(\Omega)} + \lVert Rf \rVert_{L^2(\Omega)},
$$
for some $C>0$, independent of $f$.  And lastly, $M^*$ is a continuous linear operator, so we have from  \cite{hormander2009}, Cor. 25.3.2.,

$$
\lVert f \rVert_{H^{1}(\Omega)} \le C' \lVert \chi M f \rVert_{H^{3/2}([0,T]\times S^1)} + \lVert Rf \rVert_{L^2(\Omega)}.
$$
Note that because we've multiplied $Mf(t,\theta)$ by $\chi$ and $Mf(t,\theta)$ has $\theta$ support in $S^1$, which is a compact manifold, that $\chi Mf(t,\theta)$ has compact support in $(0,T)\times S^1$, and so the norm above in $H^{3/2}([0,T]\times S^1)$ is finite.
By virtue of the injectivity of $\chi Mf(t,\theta)$, we may then write (at the cost of possibly increasing $C'$)
$$
\lVert f \rVert_{H^1(\Omega)} \le C' \lVert \chi M f \rVert_{H^{3/2}([0,T]\times S^1)}.
$$
This gives stability of the measurement operator $M$.

\begin{theorem}
Let $f\in H_D(\Omega)$, and $Mf(t,\theta)$ be defined as the either of the measurement operators above.  If $\chi \in C_0^\infty(\R)$ with $\chi \equiv 1$ on $[0,T]$, then we have the following stability estimate:
$$
\lVert f \rVert_{H^1(\Omega)} \le C\lVert \chi M f \rVert_{H^{3/2}([0,T]\times S^1)}
$$
\end{theorem}

\subsection{Visible Singularities}

A singularity $(y,\xi) \in \textnormal{WF}(f)$ is called visible on an open subset $U\times \Gamma$ of $\R \times S^1$ for $M$ if it creates a singularity in the measurement data $Mf\mid_{U\times \Gamma}$.  Now because $Mf$ is an elliptic FIO associated with a local canonical dif{}feomorphism $C$ (See Theorems \ref{thm:small},\ref{thm:large}), we know by \cite{hormander2009}
$$
\textnormal{WF}(Mf) = C \circ \textnormal{WF}(f)
$$

\begin{figure}[ht]
\centering
\includegraphics{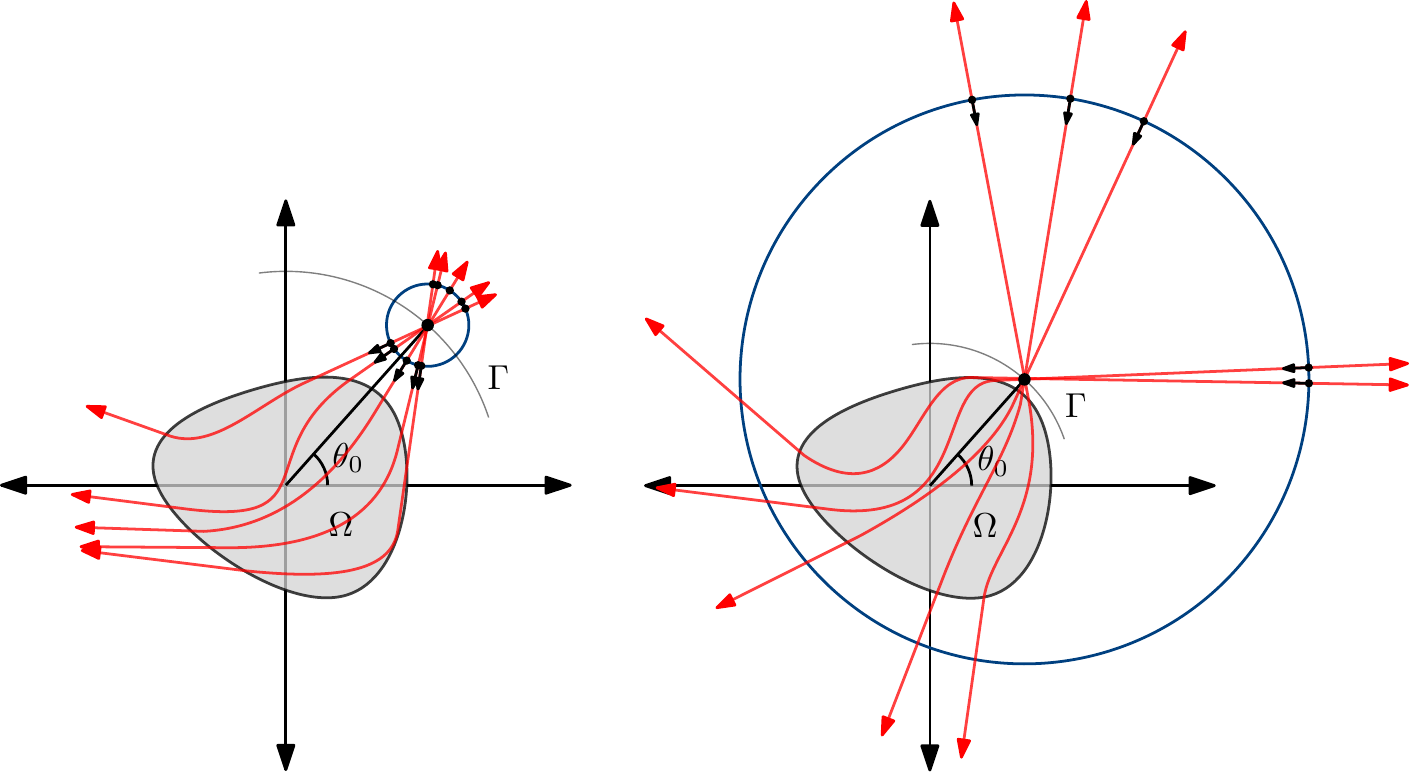} 
\caption{Singularities that may be visible from $\theta_0 \in \Gamma$ in both the cases (left) $R-r > 1$ and (right) $R=1, r>2$ will lie on the geodesics issued from the integrating detectors.}\label{fig:vis}
\end{figure}

Let $(t_0,\boldsymbol{\theta}_0) \in \R \times S^1$ and let $U\times \Gamma$ be an open neighborhood of $(t_0,\boldsymbol{\theta}_0)$.  By the above arguments (Theorems 1 and 2), we know that singularities $(y,\xi) \in \textnormal{WF}(f)$ split and travel along geodesics $\left(\gamma_{y,\hat{\xi}}(t),c(y)|\xi|\dot{\gamma}_{y,\hat{\xi}}(t)\right)$, and that this will create a singularity at $(t_0,\boldsymbol{\theta}_0)$ if and only if the geodesic intersects the circular integrating detector with center $\boldsymbol{\theta}_0$ perpendicularly at time $t_0$, with no singularity to mask it intersecting the circular integrating detector at mirror points on the circle (i.e. an antipodal point in the small radius case, and a $C(\theta,r)$ mirror point in the large radius case).  Therefore, to determine those singularities of $f$ that are visible from $U\times \Gamma$, we simply trace all geodesics that go through $R\boldsymbol{\theta}$ back to $\Omega$ and see if they have nonempty intersection with $\textnormal{WF}(f)$, see Figure \ref{fig:vis}.

For each $(t,\boldsymbol{\theta}) \in U\times \Gamma$, $\lambda \in \R\setminus 0$, let $$A^+_{t,\theta,\lambda} = \left\{ \left( \gamma_{(x,\frac{x-R\boldsymbol{\theta}}{r})}(t),\lambda \hat{\dot{\gamma}}_{(x,\frac{x-R\boldsymbol{\theta}}{r})}(t)\right) \mid |x-R\boldsymbol{\theta}| = r,\, (x-R\boldsymbol{\theta})\cdot \boldsymbol{\theta}  > 0 \right\}$$ and
$$
A^-_{t,\theta,\lambda} = \left\{ \left( \gamma_{(x,\frac{x-R\boldsymbol{\theta}}{r})}(t),\lambda \hat{\dot{\gamma}}_{(x,\frac{x-R\boldsymbol{\theta}}{r})}(t)\right) \mid |x-R\boldsymbol{\theta}| = r,\, (x-R\boldsymbol{\theta})\cdot \boldsymbol{\theta}  < 0 \right\}.$$
These are the sets of all points on geodesics intersecting the half circle $C^+(\theta,r)$ (respectively, $C^{-}(\theta,r)$) perpendicularly at time $t$, with tangent vector of magnitude $\lambda$.  For $$\left( \gamma_{(x,\frac{x-R\boldsymbol{\theta}}{r})}(t),\lambda \hat{\dot{\gamma}}_{(x,\frac{x-R\boldsymbol{\theta}}{r})}(t)\right) \in A^{\pm}_{t,\theta,\lambda},$$ define $\sim \left( \gamma_{(x,\frac{x-R\boldsymbol{\theta}}{r})}(t),\lambda \hat{\dot{\gamma}}_{(x,\frac{x-R\boldsymbol{\theta}}{r})}(t)\right) = \left( \gamma_{(\tilde{x},\frac{\tilde{x}-R\boldsymbol{\theta}}{r})}(t),\lambda \hat{\dot{\gamma}}_{(\tilde{x},\frac{\tilde{x}-R\boldsymbol{\theta}}{r})}(t)\right) \in A_{t,\theta,\lambda}^{\mp}$ where $\tilde{x}$ is the appropriate mirror point on the circular integrating detector, depending on the experimental setup.  Now, $(y,\xi) \in \textnormal{WF}(f)$ is visible from $(t,\theta) \in U \times \Gamma$ if and only if $(y,\xi) \in A^{\pm}_{(t,\theta,\lambda)}$ for some $(t,\theta,\lambda) \in U\times \Gamma \times (\R \setminus 0)$ and $\sim (y,\xi) \not \in \textnormal{WF}(f)$.  Let 
$$B^{\pm}_{t,\theta,\lambda} = \left\{ (y,\xi) \mid (y,\xi) \in A^{\pm}_{t,\theta,\lambda} \cap \textnormal{WF}(f) \text{ and } \sim (y,\xi) \not\in \textnormal{WF}(f)\right\}.$$
It then follows from the above arguments that the set of visible singularities is given by
$$
\bigcup_{(t,\theta,\lambda) \in U \times \Gamma \times (\R \setminus 0)} B^+_{t,\theta,\lambda} \cup B^{-}_{t,\theta,\lambda}.$$
We have shown the following:
\begin{theorem}
Let $U \times \Gamma \subset \R \times S^1$ be an open subset, and for each $(t,\theta,\lambda) \in U\times \Gamma\times (\R \setminus 0)$ let $A^{\pm}_{t,\theta,\lambda}$ and $B^{\pm}_{t,\theta,\lambda}$ be defined as above.  Then in both the small radius detector case and the large radius detector case, the singularities of $f$ that are visible from $U\times \Gamma$ in the restricted data $Mf\mid_{U\times \Gamma}$ are given by 
$$
\bigcup_{(t,\theta,\lambda) \in U\times \Gamma\times (\R \setminus 0)} B_{t,\theta,\lambda}^+ \cup B_{t,\theta,\lambda}^-.
$$
\end{theorem}

From this we see that if $\Gamma = S^1$ and $(0,T] \subset U$, where $T = \sup_{x\in\Omega} \inf_{\theta\in[0,2\pi)} \text{dist}(x,C(\theta,r))$ (where the distance is the geodesic distance), then all singularities of $f$ are visible assuming $T < \infty$.  We make a note that just because a singularity is not ``visible'' from some open set $U \times \Gamma$ does not mean that we cannot infer the existence of that singularity from the restricted data $Mf$.
\section{Numerical Results}
\begin{figure}[h!]
\centering
\includegraphics[scale=0.72]{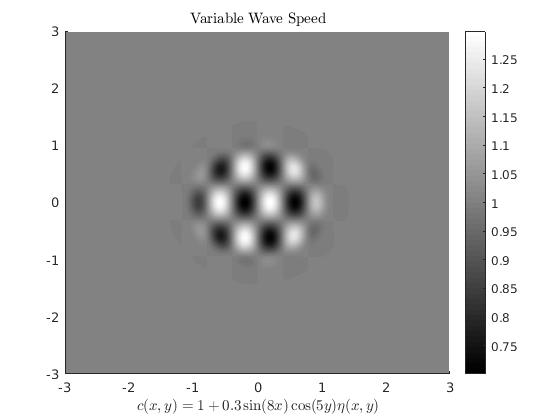} 
\caption{Variable wave speed of $1+0.3\sin(8x)\cos(5y)\eta(x,y)$, where $\eta(x,y)\in C_0^\infty(B_1(0))$.}\label{fig:waveSpeed}
\end{figure}

\begin{figure}[h!]
 \centering
 \includegraphics{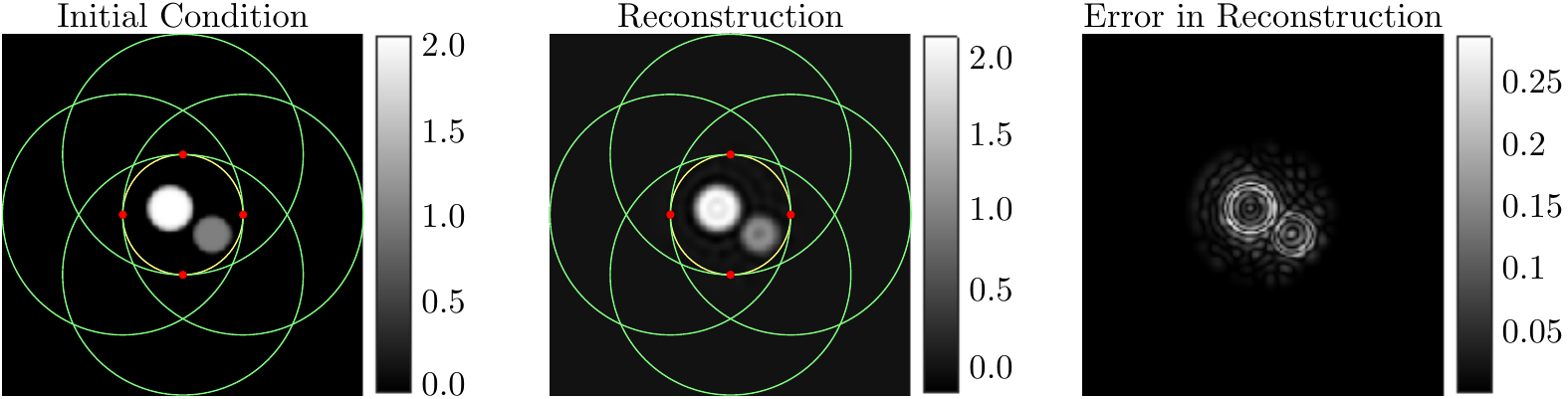}
 \caption{Results of reconstruction using $R=1$ and $r=2$ model (Large radius detector model).  This reconstruction was made using full data.}\label{fig:discs}
\end{figure}

To simulate the collection of forward data, we numerically solve the wave equation with variable wave speed using the implementation of Perfectly Matched Layers (PML) found in \cite{grote} for a number of dif{}ferent smooth initial conditions.  This ensures that measured data will only come from signals inside the region of interest, and not from reflections at the boundary of the window of computation.  Then, we collect simulated measurement data on the unit circle $Mf(t,\theta)$ for $0 \le \theta < 2\pi$ and $0 \le t < 5$, for a specific initial condition.  In general, the amount of time that we collect data should depend on the wave speed inside the medium we are imaging, and for the wave speed we have chosen of $1+0.3\sin(8x)\cos(5y)\eta(x,y)$ with $\eta(x,y)\in C_0^\infty(B_1(0))$, t=5s suf{}fices as an appropriate time range.  We've show the graph of the wave speed in Figure \ref{fig:waveSpeed}.  We then use an iterative solver to reconstruct the smooth initial condition using the simulated data over the given time interval.  The reconstruction shown in Figures \ref{fig:discs} was made using the model $R=1$ and $r=2$ (the large radius detector model), with data taken on the full unit circle.  An almost identical reconstruction is obtained if we use the small radius integrating detector model with full data.

\begin{figure}[h!]
    \centering
    \includegraphics{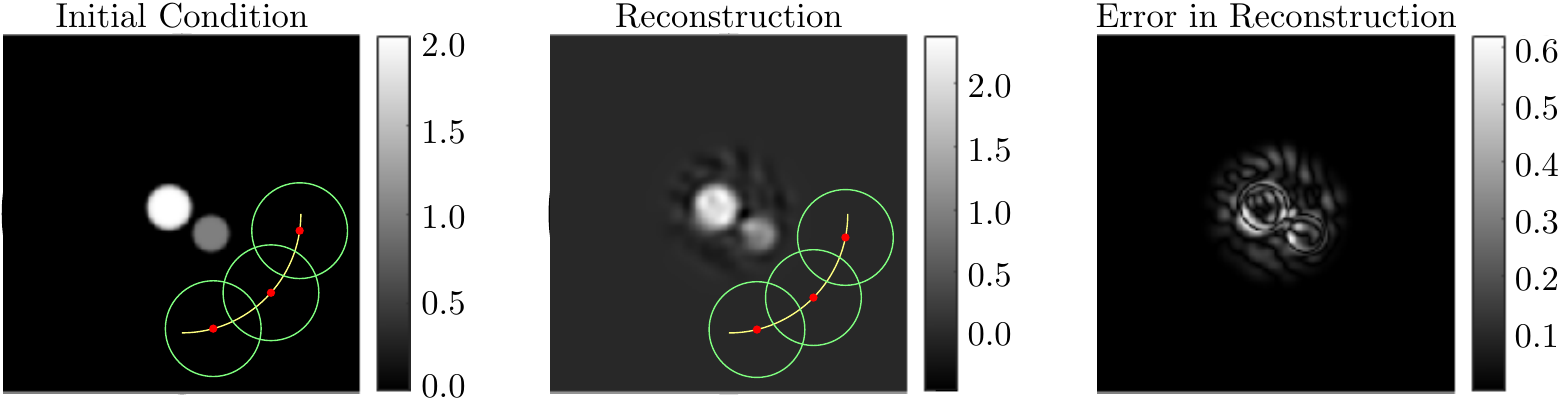}
    \caption{Result of reconstruction with partial data using $R = 2$, and $r=0.8$ (Small radius detector model).  This reconstruction was for $\theta \in (-\pi/2, 0)$. Shown in the figure are the set on which data is collected as well as some representative circular integrating detectors.}
    \label{fig:part_S}
\end{figure}

We also run numerical simulations with data taken on an open subset of the unit circle.  Here we take data for $\theta \in (-\pi/2,0)$, with the same wave speed interior to the object.  We then multiply by a smooth cutof{}f function so as to not introduce new singularities into the reconstruction.  Results for the partial data case are shown in Figure \ref{fig:part_S}.

\begin{figure}[h!]
    \centering
    \includegraphics{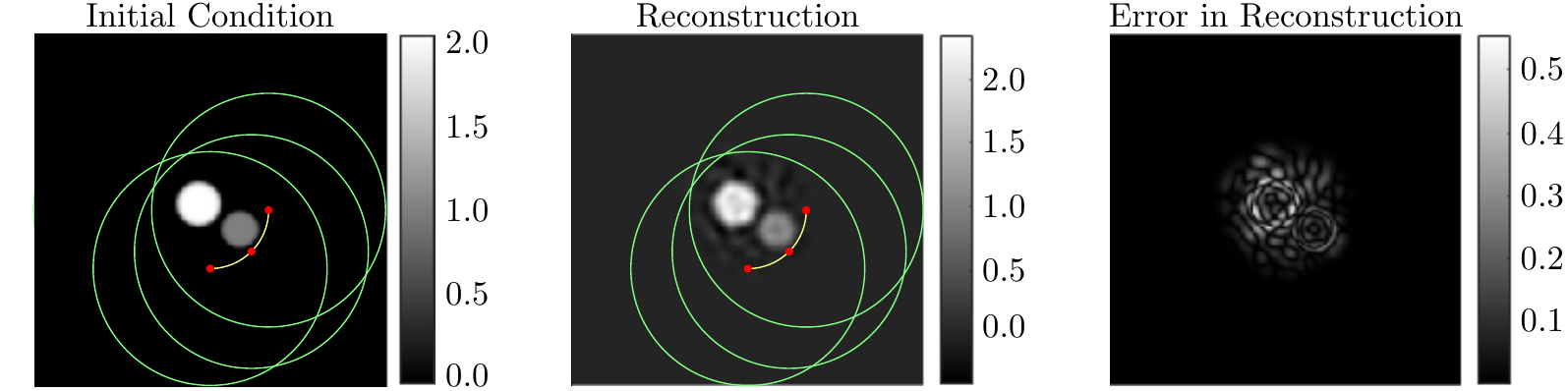}
    \caption{Result of reconstruction with partial data using $R=1$, and $r=2$ (Large radius detector model).  This reconstruction was for $\theta \in (-\pi/2,0)$. Shown in the figure are the set on which data is collected as well as some representative circular integrating detectors.}
    \label{fig:part_B}
\end{figure}

\printbibliography
\end{document}